# On The Ten Point Wood-Desargues' Configuration

**Christopher J Bradley**

## Abstract

The ten pairs of directly similar triangles in perspective in the Wood-Desargues' configuration have ten common Hagge circle centres. These centres are known to lie one each of the ten perspectrices. It is shown in this paper that these centres also lie at the vertices of five cyclic quadrangles, which are similar to the five cyclic quadrangles of the initial Wood-Desargues' configuration and that these five new cyclic quadrangles have a common circumradius.

## 1. Introduction

The Wood-Desargues' configuration is built up as follows: Start with two circles intersecting at points *J* and *K*. In the first circle inscribe a triangle *ABC*. Draw *AK*, *BK*, *CK* to meet the second circle at points *a*, *b*, *c*. It follows that triangles *ABC* and *abc* are directly similar triangles in perspective with *J* the centre of direct similarity. By this we mean that triangle *abc* may be obtained by from triangle *ABC* by a rotation about *J* followed by a dilation centre *J*. Suppose that *BC* meets *bc* at the point 1, *CA* meets *ca* at the point 2 and *AB* meets *ab* at the point 3. Then 123 is the Desargues' axis of perspective, commonly referred to nowadays as the perspectrix. There are in all ten points *A*, *B*, *C*, *K*, *a*, *b*, *c*, 1, 2, 3 which form what may be called the Wood-Desargues' configuration. The ten points turn out to be equivalent in the following sense. Each of them is a vertex of perspective of two triangles in the configuration, the remaining three points being collinear and serving as the perspectrix. This situation is illustrated in Fig. 1. In his landmark paper Wood [1] proved that the quadrangles *Aa*23, *Bb*31, *Cc*12 are cyclic in addition to the starting quadrangles *ABCK* and *abcK*. These circles are shown in Fig. 1 and for reasons, which will emerge shortly, their centres are also labelled *ABCK*, *abcK* etc. In his paper Wood showed that the five circle centres are themselves concyclic and lie on a circle passing through *J*. The ten pairs of triangles in perspective, together with their vertices of perspective and their perspectrices are shown in Table 1. The orders of the vertices in pairs of triangles correspond and the order of letters in the perspectrix is standard, so, for example, 1 is the intersection of *BC* and *bc*, *BC* and 32 and *bc* and 32.

| Triangle 1 | Triangle 2 | Vertex | Perspectrix |
|---|---|---|---|
| *ABC* | *abc* | *K* | 123 |
| *KBC* | *a*32 | *A* | 1*cb* |
| *AKC* | 3*b*1 | *B* | *c*2*a* |
| *ABK* | 21*c* | *C* | *ba*3 |
| *Cc*2 | *Bb*3 | 1 | *aAK* |
| *Aa*3 | *Cc*1 | 2 | *bBK* |
| *Bb*1 | *Aa*2 | 3 | *cCK* |
| *Kbc* | *A*32 | *a* | 1*CB* |
| *Kca* | *B*13 | *b* | 2*AC* |
| *Kab* | *C*21 | *c* | 3*BA* |

**Table 1: The ten perspectives**

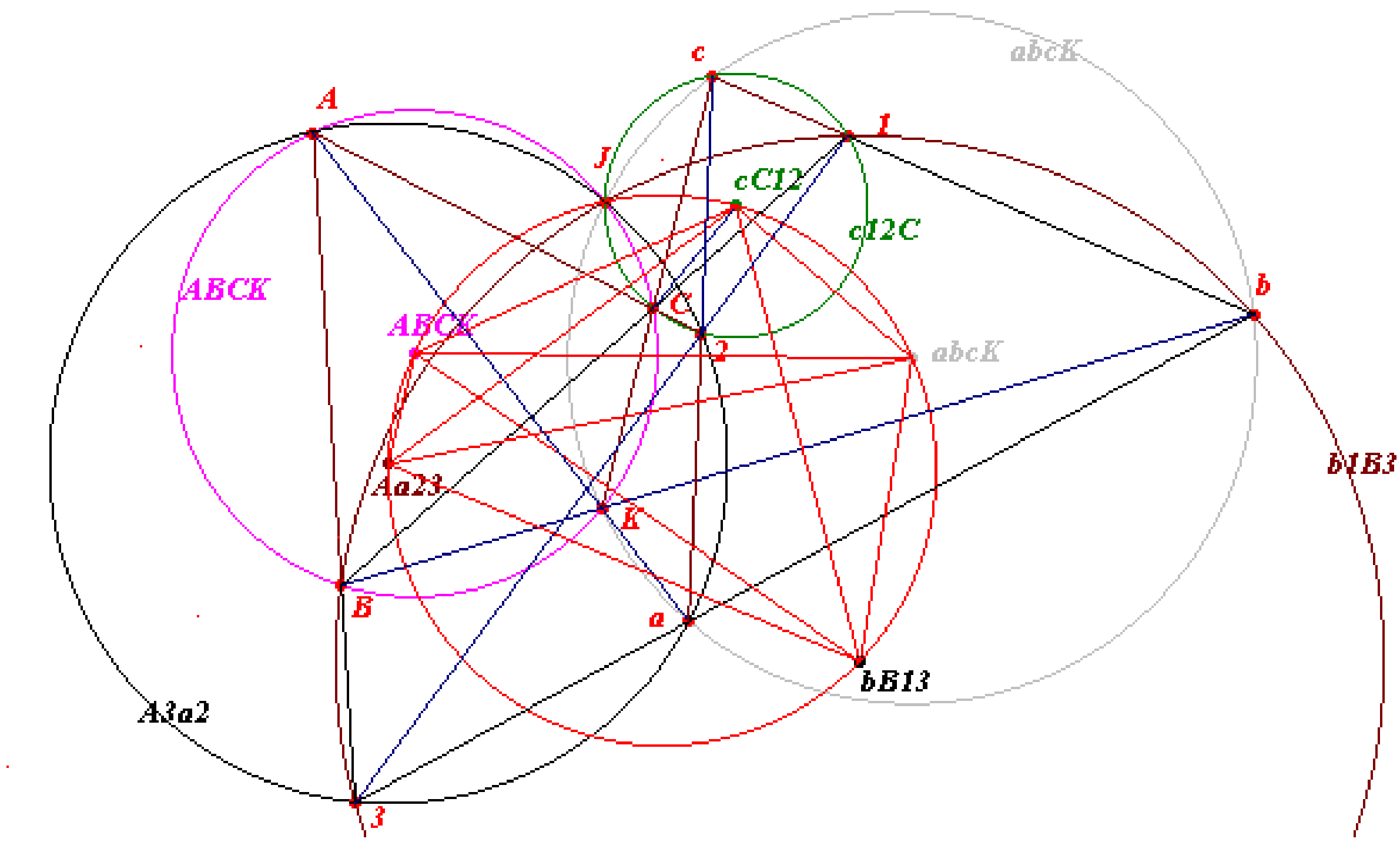

**Fig.1**

**The Wood Desargues' configuration**

In Section 2 we consider the orthocentres of the twenty triangles and show that they lie four at a time on five quadrangles similar to the initial five quadrangles and also that the same twenty orthocentres lie four at a time on five lines. This follows from well-known results about cyclic quadrangles and quadrilaterals, but is presented, as the result is required in what follows.

In Section 3 we establish our main results. Each pair of triangles in the list of ten perspective pairs have orthocentres, say *H* and *F*, so that the circle *JHF* is a common Hagge circle of the two triangles. Its centre we denote by *h*, and it is the intersection of the perpendicular bisectors of *JH* and *JF*. It is known that these ten points *h* lie one each on the ten perspectrices, see [1]. We show that they also form the vertices of five cyclic quadrangles similar to the initial five cyclic quadrangles, and that each of the five circles involved has the same circumradius. The method of proof involves the pentagon formed by the five circle centres of the initial five circles. The result follows from lemmas that are proved and from straightforward extensions of the work of Wood [1].

## 2. The twenty orthocentres

In Fig. 2 we show the circle *ABCK* and its four component triangles *ABC*, *KBC*, *AKC*, *ABK* with orthocentres *H*(*K*), *H*(*A*), *H*(*B*), *H*(*C*) respectively. In the same figure we also show the four lines 123, *b*1*c*, *ac*2, 3*ba*, which are the four perspectrices of the four perspectives associated with the triangles of the cyclic quadrangle *ABCK*. The four orthocentres of the triangles formed by these four lines, namely triangles *abc*, *a*32, 3*b*1 and 21*c* (see Table 1) are denoted by *F*(*K*), *F*(*A*), *F*(*B*), *F*(*C*) respectively.

*Proposition* 1
The quadrangle *H*(*A*)*H*(*B*)*H*(*K*)*H*(*C*) is congruent to the quadrangle *ABKC*.

*Proof*
The result is true for any cyclic quadrangle and does not depend upon *A*, *B*, *C*, *K* being part of the Wood-Desargues' configuration. Let the centre of the circle *ABCK* be the origin of vectors and suppose *A*, *B*, *C*, *K* have vector positions $\boldsymbol{a}$, $\boldsymbol{b}$, $\boldsymbol{c}$, $\boldsymbol{k}$ respectively. Then *H*(*K*) has vector position $\boldsymbol{a} + \boldsymbol{b} + \boldsymbol{c}$ and *H*(*A*) has vector position $\boldsymbol{b} + \boldsymbol{c} + \boldsymbol{k}$. It follows that $\boldsymbol{H}(\boldsymbol{A})\boldsymbol{H}(\boldsymbol{K}) = \boldsymbol{a} - \boldsymbol{k} = \boldsymbol{KA}$. It is evident that *H*(*A*)*H*(*B*)*H*(*K*)*H*(*C*) is the image of *ABKC* under a 180$^{\circ}$ rotation about the point whose vector position is ½ ($\boldsymbol{a} + \boldsymbol{b} + \boldsymbol{c} + \boldsymbol{k}$).

*Proposition* 2
The points *F*(*A*), *F*(*B*), *F*(*C*), *F*(*K*) are collinear.

*Proof*
The result is true for the orthocentres of the four triangles formed by any set of mutually non-parallel lines and does not depend on the lines 123, *b*1*c*, *ac*2, 3*ba* being part of the Wood-Desargues' configuration. The result is immediate from the work of Steiner on the complete quadrilateral. See, for example, Durrell [2].

Now there are five circles and associated quadrangles in the initial Wood-Desargues' configuration, so there are similar results to Propositions 1 and 2 for each of the other four. This means that each of the twenty orthocentres of the twenty triangles of the ten perspectives appears twice. Each will appear once as an *H* point of one of the five quadrangles and once as an *F* point of one of the five lines. The five resulting cyclic quadrangles are congruent in each case to one of the initial five quadrangles.

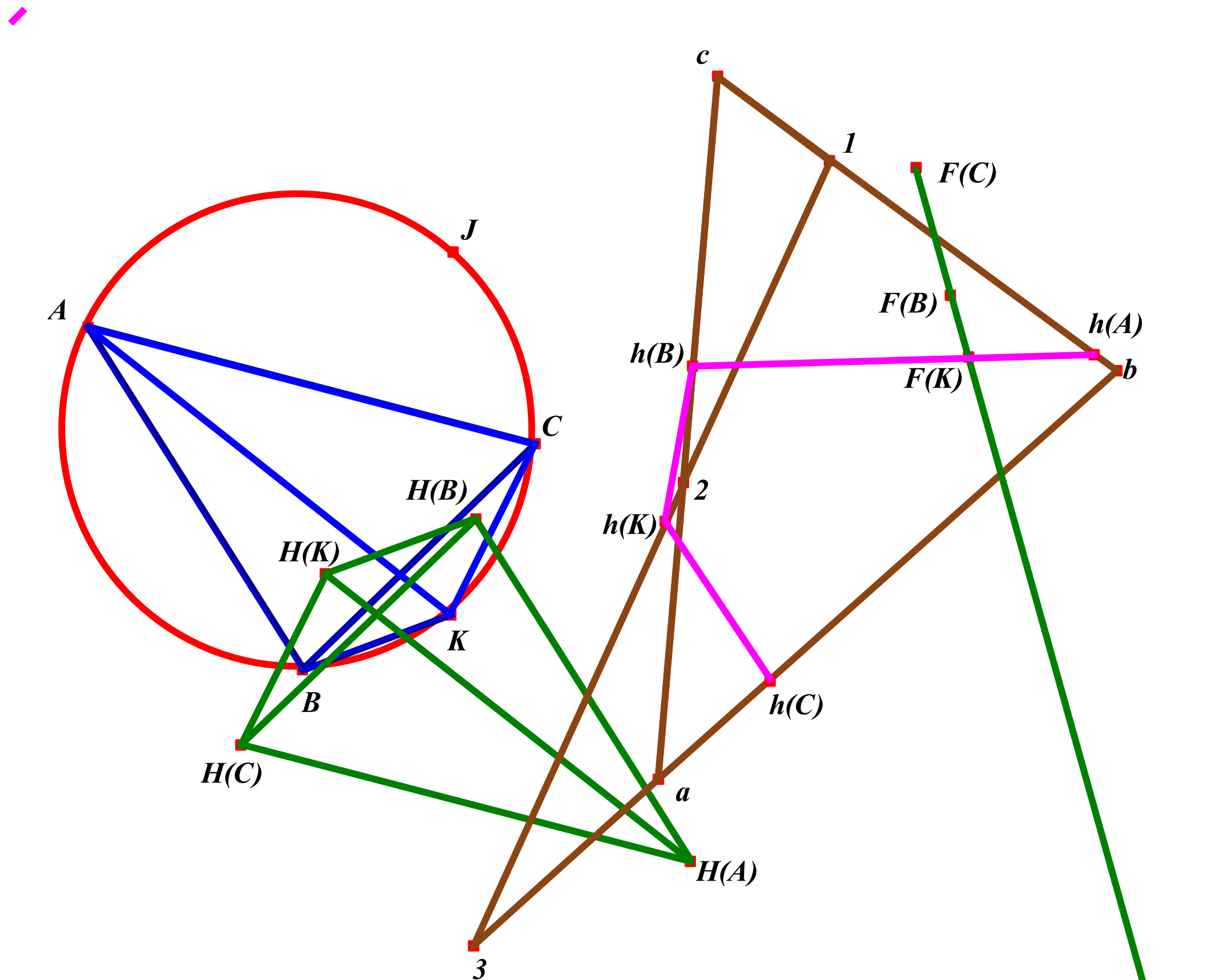

**Fig. 2**
**ABCK, its triangles, perspectives, perspectrices and orthocentres**

The perpendicular bisectors of $JH(K)$ and $JF(K)$ meet at a point we denote by $h(K)$. This is the joint Hagge centre of triangles $ABC$ and $abc$. It is known to lie on the perspectrix 123. The points $h(A)$ and $h(B)$ and $h(C)$ are similarly defined for the pairs of triangles $KBC$, $a32$ and $AKC$, $3b1$ and $ABK$, $21c$ respectively and lie on the lines $b1c$, $ac2$, $3ba$ respectively. In Fig. 2 we also show these points and the quadrangle $h(A)h(B)h(K)h(C)$. There are five such quadrangles and each of the $h$ points appears on two of them. For example $h(K)$ lies on the quadrangle $h(a)h(b)h(c)h(K)$ as well. In Section 3 we show that these five quadrangles are directly similar to the five initial quadrangles in the Wood-Desargues' configuration.

## 3. The ten common Hagge centres

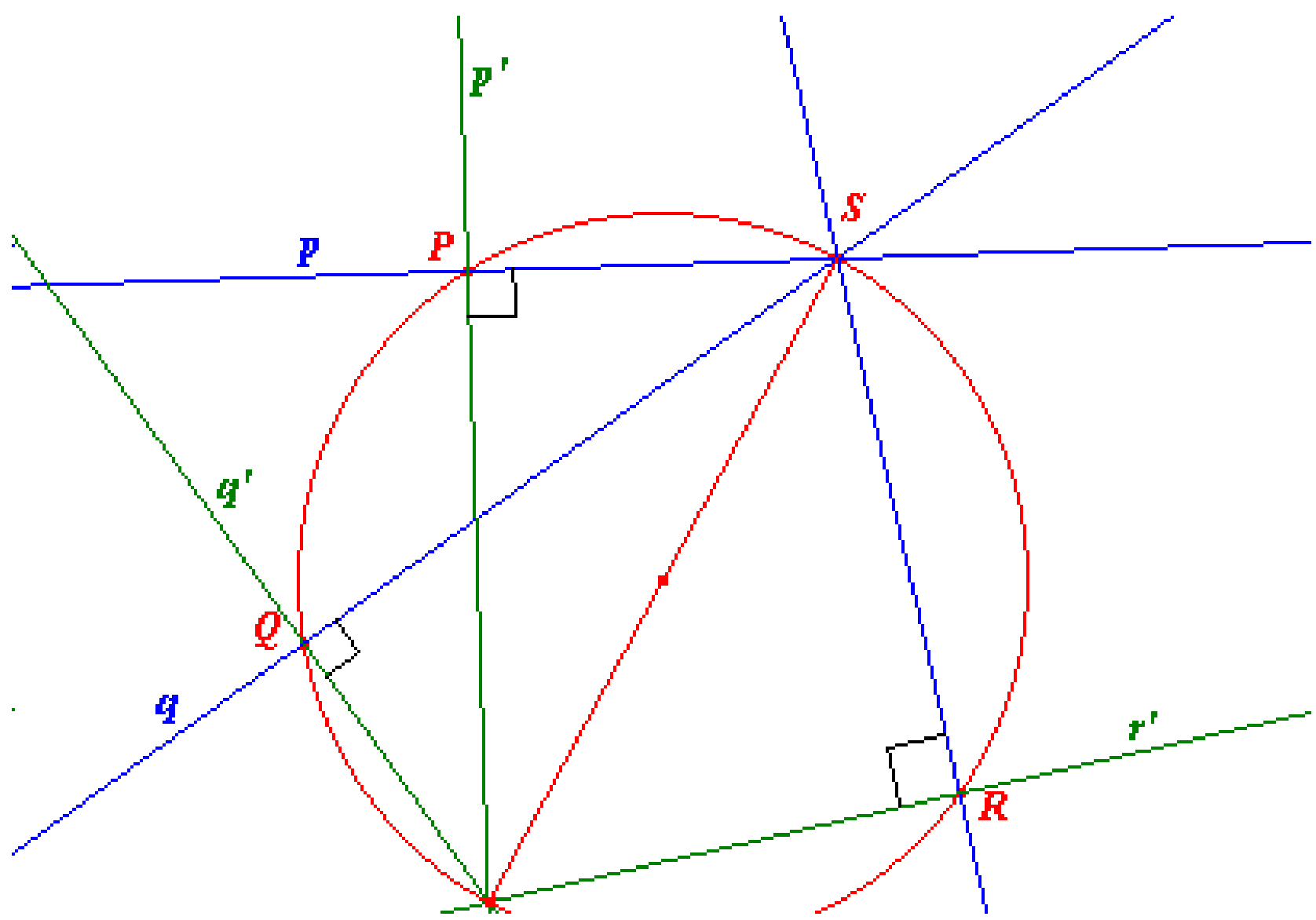

**Fig. 3**
**Illustration of Lemma 1**

First we establish a very simple Lemma that at first sight does not appear to be related the Wood-Desargues' configuration at all, but is in fact a useful step in establishing two of the theorems we want to prove. The Lemma is illustrated in Figure 3.

*Lemma* 1
Suppose that $P$, $Q$, $R$ are three non-collinear points and lines $p$, $q$, $r$ through $P$, $Q$, $R$ respectively meet at a point $S$ on circle $PQR$. Then lines $p'$, $q'$, $r'$ through $P$, $Q$, $R$ respectively, perpendicular to $p$, $q$, $r$ respectively, are concurrent at a point $T$ diametrically opposite to $S$.
*Proof*
Let the lines $p'$ and $q'$ meet at $T$. Then, since $\angle SPT = \angle SQT = 90^{\circ}$, $T$ lies at the other end of the diameter to $S$ of circle $SPQ$. But $R$ lies on this circle, and $r$ and $r'$ are at right angles. It follows that $r'$ passes through $T$.

It is not needed here, but it is easy to verify that the condition that $S$ lies on circle $PQR$ is necessary as well sufficient for the concurrence of $p'$, $q'$, $r'$.

We now establish the crucial Lemma in helping to establish what we wish to prove. One possible configuration is shown in Fig. 4, but the proof given using directed angles is diagram independent.

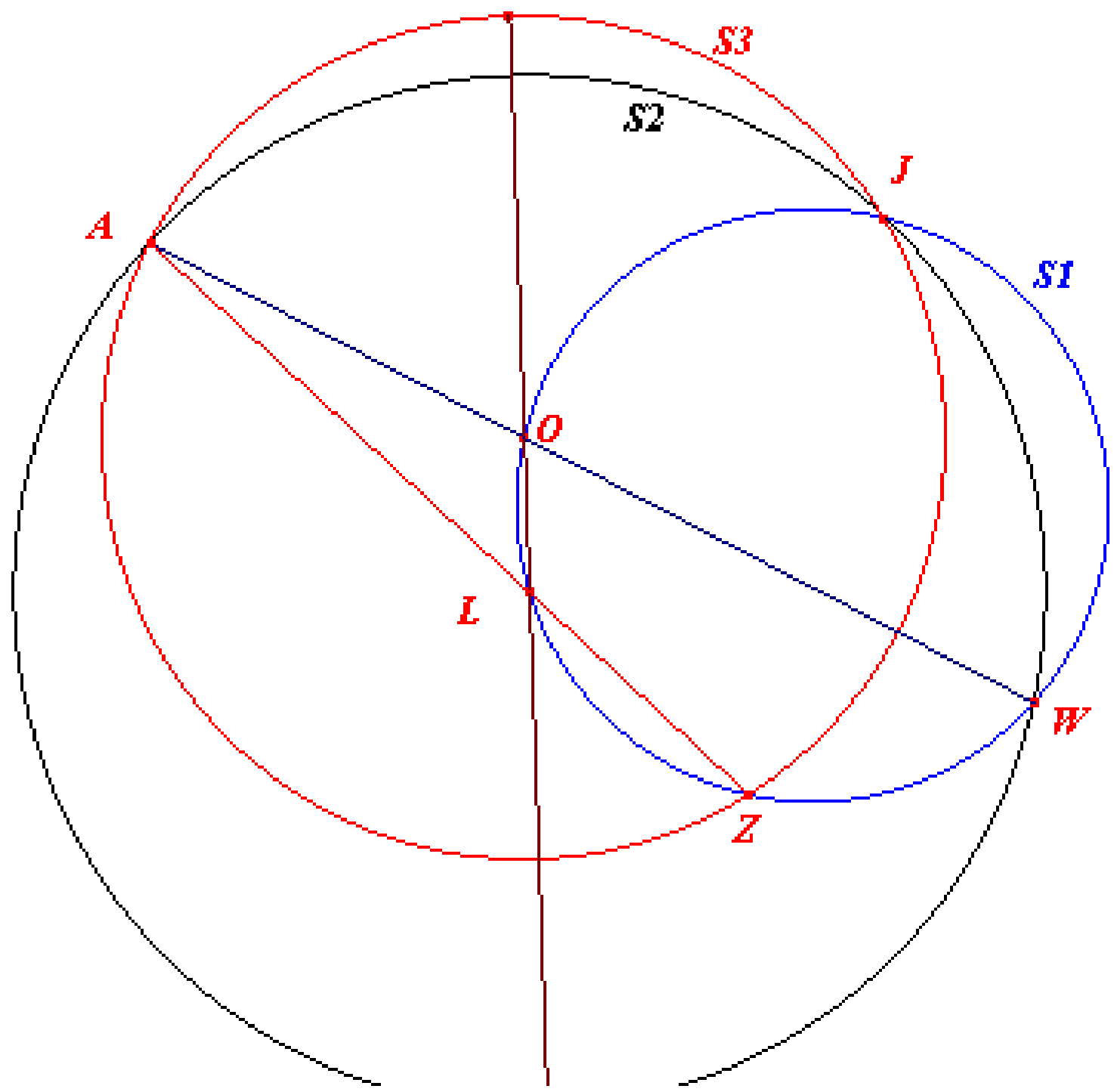

**Fig. 4**
**Illustration of Lemma 2**

*Lemma* 2
Let *J*, *O*, *L* be three points on the circumference of a circle $S_1$ and let $S_2$ and $S_3$ be the circles through *J*, centres *L* and *O* respectively. Suppose that $S_2$ and $S_3$ meet again at *A*, $S_3$ and $S_1$ meet again at *Z* and $S_1$ and $S_2$ meet again at *W*, then *A*, *L* and *Z* are collinear as are *A*, *O* and *W*.

*Proof*
Fig. 5 is a copy of Fig. 1, but with some additional points and lines added. First, however, note that the centres of the five circles have been relabelled, so that circle *ABCK* has centre *U*, circle *abcK* has centre *V*, circle *Aa23* has centre *L*, circle *Bb13* has centre *M* and circle *cC12* has centre *N*. Consider the point *L*, the centre of circle *Aa*23. This is clearly the point of intersection of the perpendicular bisectors of *J*A, *J*a, *J*2, *J*3.

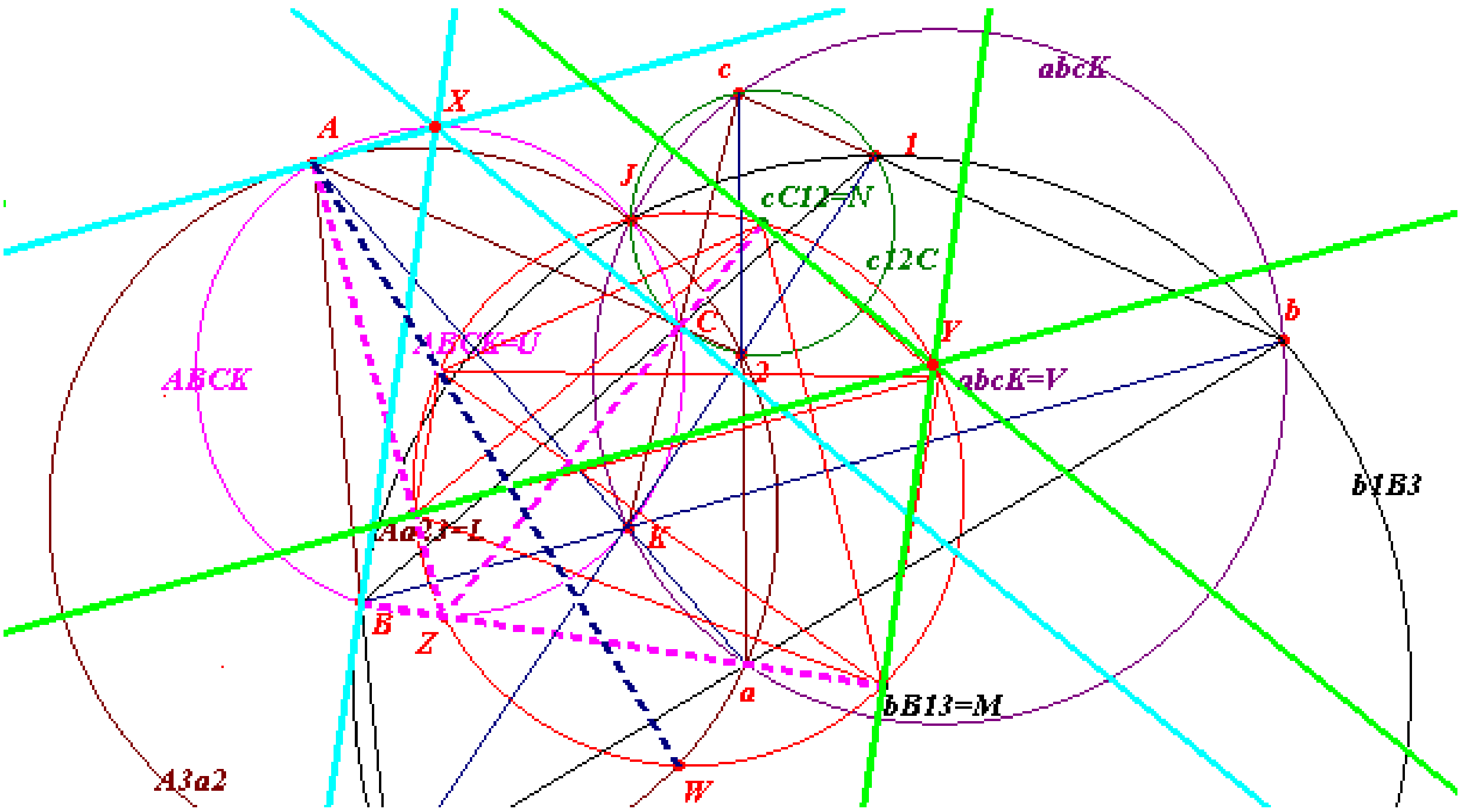

**Fig. 5**
**Fig. 1 with extra properties illustrated**

Now *A, a*; *A*, 2; *A*, 3; *a*, 2; *a*, 3, 2, 3 are six corresponding pairs of vertices amongst the thirty such pairs in the ten perspectives in Table 1. The same applies to each of the five circle centres and this fact accounts for all thirty such points of intersection.

Additional points in the diagram are *X*, *Y*, *Z*, *W*. These points occur in Propositions we establish next. Lines *ALZ*, *AUW*, *BMZ*, *CNZ* are also shown, and these collinearities will soon be established.

*Proposition* 3
Let *Z* be the point of intersection of circles *JUVLMN* and *JABCK*, and let *W* be the point of intersection of circles *JUVLMN* and *JAa*23, then *ALZ*, *AUW*, *BMZ*, *CNZ* are straight lines.

*Proof*
In any ten point Wood-Desargues' configuration, the conditions of Lemma 2 apply, which accounts for the straight lines *ALZ* and *AUW*. Application of the same Lemma to other pairs of circles shows that *BMZ* and *CNZ* are also straight lines.

*Proposition* 4
Triangles *ABC* and *LMN* are directly similar triangles in perspective with centre of direct similarity *J* and vertex of perspective *Z*.

*Proof*
This follows immediately as a result of Proposition 3.

Since circle *JUVLMN* cuts the five initial circles in five different points, there are thus five such vertices of perspective and all twenty triangles of the initial Wood-Desargues' configuration get mapped by perspectives on to the triangles of the pentagon *UVLMN*. We have thus established

*Proposition* 5
Consider the pentagon *UVLMN*. Omit one vertex, say *M*. Then the quadrangle *UVLN* is directly similar to the quadrangle *Bb*31 (a little care has to be taken over the order of the vertices). The same applies when any other vertex of the pentagon is omitted.

In this way the pentagon inscribed in the circle *UVLMN* carries all the angular information of the five quadrangles of the initial configuration.

If we now use Lemma 1 with *X* the point diametrically opposite *Z* on the circle *JABCK* we recover Theorem 12 of Wood [1], namely

*Proposition* 6
Tangents at *A*, *B*, *C* to circles *Aa*23, *Bb*13, *Cc*12 are concurrent at a point *X* on circle *JABCK*.

*Proof*
These tangents are perpendicular to the normals *AL*, *BM*, *CN* so Lemma 1 applies.

There are, of course other cases in which sets of tangents are concurrent.

Lemma 1 also provides proof of the following proposition.

*Proposition* 7
Lines through *L*, *M*, *N* parallel to the tangents at *A*, *B*, *C* in Proposition 6 are concurrent at a point *Y* on circle *JUVLMN*, which is diametrically opposite *Z* in that circle.

This allows us to identify the centre of the circle *JUVLMN*. It is the midpoint of *YZ*.

Recall now that *h* is the intersection of the perpendicular bisectors of *JH* and *JF* and this is how all ten common Hagge centres are constructed. We have shown above that if we take the perpendicular bisectors of *JA* and *Ja* we obtain *L*, if we take the perpendicular bisectors of *JB* and *Jb* we get *M* and if we take the perpendicular bisectors of *JC* and *Jc* we get *N*. But we have shown in Proposition 5 that triangles *ABC* and *LMN* are in perspective and similar with vertex. It follows that *h*(*K*) is the orthocentre of triangle *LMN*. Similar considerations hold for all the other common Hagge centres, which together account for all ten orthocentres of the triangle formed by the pentagon *UVLMN* and its diagonals.

*Proposition* 8
The quadrangle $h(A)h(B)h(C)h(K)$ is similar to the quadrangle $ABCK$.

*Proof*
We know from Proposition 1 that the quadrangles $ABCK$ and $H(A)H(B)H(C)H(K)$ are congruent. For the same reason quadrangles $LMNV$ and $h(A)h(B)h(C)h(K)$ are congruent. But from Proposition 5 the quadrangles $ABCK$ and $LMNV$ are similar.

There are five such cyclic quadrangles formed by the common Hagge centres and since the orthocentres of the triangles formed by $UVLMN$ and its diagonals all derive from the same circle we have proved

*Proposition* 9

The ten common Hagge centres of the Wood-Desargues' configuration are inscribed in five circles of equal radius and each Hagge centre lies on two such cyclic quadrangles.

The configuration is illustrated in Fig. 6, in which the initial configuration and its five circles and quadrangles are shown, together with the perspectrices and the common Hagge centres lying on the perspectrices and forming five cyclic quadrangles of equal radius.

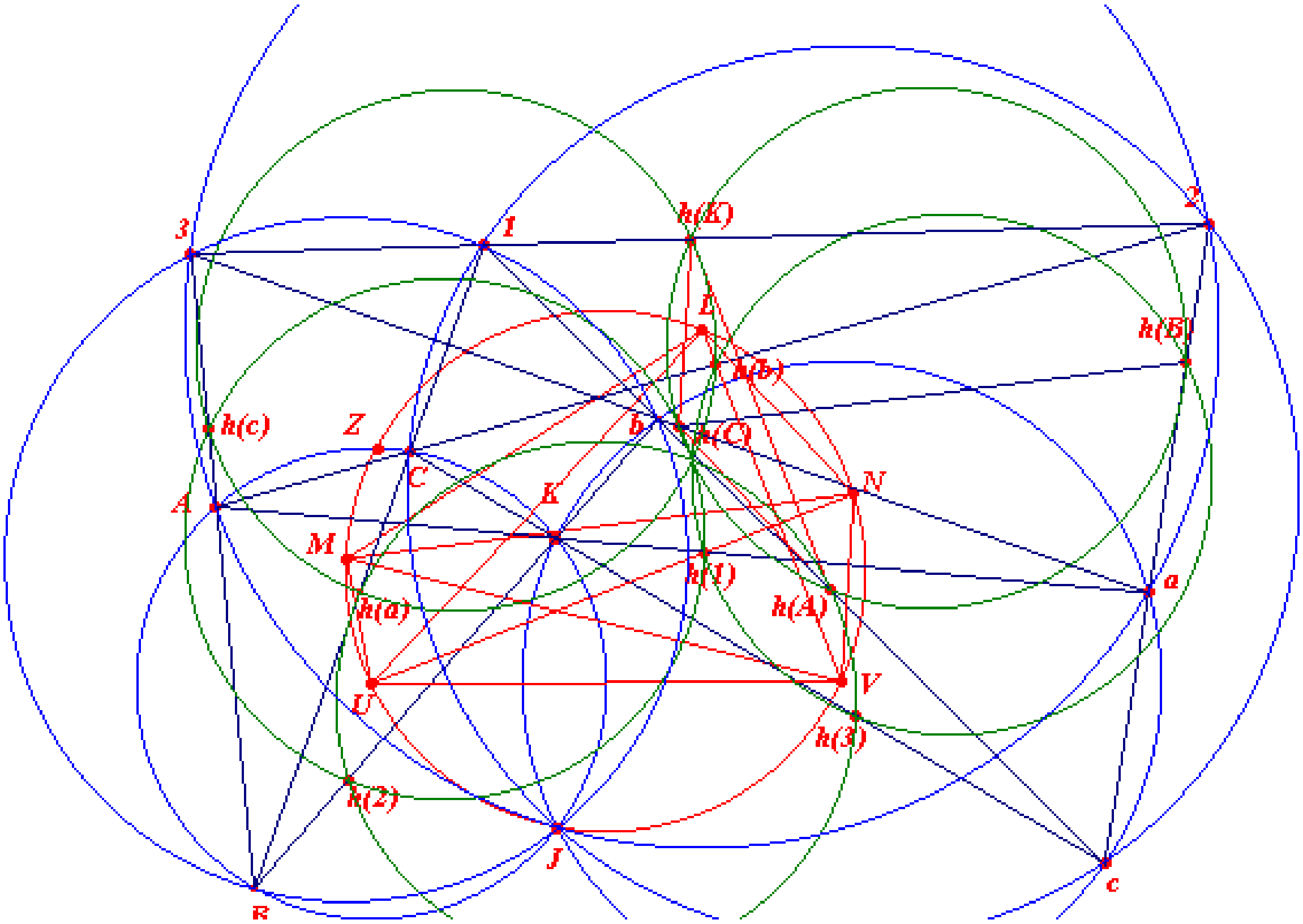

**Fig. 6**
**The five circles and quadrangles**

*Acknowledgement*
Thanks are due to Geoff. C. Smith (of the University of Bath) who helped me substantially with the work in this article. If sent to a journal it would have been a joint publication.